\documentclass[12pt]{article}

\RequirePackage{a4}
\usepackage{epsfig}
\usepackage[utf8x]{inputenc}
\usepackage{verbatim}
\usepackage{hyperref}
\usepackage{url}

\begin{document}

\title{On the Japanese Multiplication Method\footnote{Published in
   {\em Progetto Alice,
     rivista di matematica e didattica},
   nr. 55 (2018) pp. 55-72.} \\
  --\, A father-and-daughter dialogue \,--  \\
(In Italian)}

\author{G.~D'Agostini \\
Universit\`a ``La Sapienza'' and INFN, Roma, Italia \\
{\small (giulio.dagostini@roma1.infn.it,
 \url{http://www.roma1.infn.it/~dagos})}
}

\date{}

\maketitle

\begin{abstract}
Recently the media broadcast the news, together with illustrative
videos, of a so-called Japanese method to perform multiplication by
hand without using the multiplication tables.\footnote{A Google
  search of the key words ``japanese method multiplication''
produced about 80 millions results at the beginning of February 2018.}
``Goodbye multiplication
tables'' was the headline of several websites, including important
ones, where news are however too often `re-posted' uncritically. The
easy numerical examples could induce naive internauts to believe that,
in a short future, multiplications could be really done without the
knowledge of multiplication tables. This is what a girl expresses,
with great enthusiasm, to her father. The dialogues described here, although
not real, are likely and have been inspired by this episode, being
Maddalena the daughter of the author. Obviously the revolutionary
value of the new method is easily disassembled, while its educational
utility is highlighted to show (or remember) the reasoning on which
the method learned in elementary school is based, although mostly
applied mechanically.
\end{abstract}

\newpage
\begin{center}
  {\Huge Il metodo giapponese delle moltiplicazioni } \\
  \mbox{} \\
       {\Huge --\,  Un dialogo fra padre e figlia\, --}
       \mbox{} \\  \mbox{} \\  \mbox{} \\ \mbox{} \\
            {\bf Sommario} 
\end{center}
Recentemente i media hanno passato la notizia, accompagnata da video
illustrativi,\footnote{Il link \url{https://www.facebook.com/didyouknowpage1/videos/1549535068487238/}
  al video citato nella prima versione di questo dialogo (inizio 2018)
  è risultato
    successivamente {\em broken}. Ma ci sono altri video su youtube
    a cui si arriva
    cercando, ad esempio, le parole chiave `metodo giapponese moltiplicazioni'. }
di un cosiddetto {\em metodo giapponese} per eseguire
moltiplicazioni a mano senza far uso delle tabelline.\footnote{Una
  ricerca su Google 
  delle parole chiave ``japanese method multiplication''
  produceva circa 80 milioni di risultati all'inizio di febbraio 2018
  (``metodo giapponese moltiplicazioni'', in italiano, ne dava
  {\it solo} circa 32 mila). Il 3 dicembre 2022, data di rivisitazione
  del testo, ``japanese method multiplication'' dava 2.110.000 risultati,
  ``metodo giapponese moltiplicazioni'' 6980.
 }
``Addio tabelline'' titolavano addirittura siti vari, anche di testate
influenti\,\cite{Focus},
nei quali le notizie vengono troppo spesso `ripostate' in
modo acritico.
 I facili esempi numerici potevano indurre ingenui
internauti a credere che in un prossimo futuro si potesse fare a meno
delle tabelline per eseguire moltiplicazioni a mano. È quanto
manifesta una ragazza con grande entusiasmo al padre, episodio dal
quale prendono spunto il dialogo qui riportato, il quale pur non
essendo reale, è verosimile, essendo Maddalena la figlia
dell'autore. Ovviamente la valenza rivoluzionaria del nuovo metodo
viene facilmente smontata, mentre ne viene messa in luce la sua
utilità didattica per mostrare, o ricordare, i ragionamenti sui quali
si basa il metodo appreso alle elementari ed eseguito per lo più
meccanicamente.
  
\section*{}
\begin{description}
 \setlength\itemsep{0.5mm}
\item[F.] Guarda che figo, papà! 
\item[P.] Qualche impresa da Guinness? O nuovo modello di\ldots
\item[F.] No, no, stavolta è roba seria, e ti interesserà di sicuro!
\item[P.] Una novità che possa interessare sia a te che a me?
  È proprio il caso di dire
  ``o fatto proprio strano da giocarci al lotto\ldots''\footnote{
    Primi versi de {\em Il cappotto}, di Gianni Rodari \cite{Rodari}.}
   Non succedeva da quando usciva un nuovo gioco della Wii!
\item[F.] In effetti è quasi un gioco.
\item[P.] Aspetta un istante che invio questo mail e\ldots
\item[F.] {\em Questa} mail, papà. Mail è femminile. Lo ha detto pure\ldots 
\item[P.] Senti, Maddale', io mando mail da più di trent'anni e me ne frego
  di quello che ha decretato qualche sapientone.
\item[F.] Ma `posta' è femminile\ldots
\item[F.] E chi è che manda `una posta'? Mandiamo un messaggio,  e forse
  è per questo che nella mia mente `mail' è maschile. E infatti, nel mio ambiente
  si diceva ``an e-mail message''. Comunque, ecco, un momento,\ldots\ 
  {\em e-mail sent}, così
  hai poco da ridire -- e, ora che mi ci fai pensare, `mail'
  in inglese è neutro e noi rendiamo i nomi neutri con il
  maschile\ldots \ ma lasciamo stare\ldots \footnote{Il lettore
    potrebbe essere interessato a cosa ne pensa
    l'Accademia della Crusca\,\cite{Crusca} (nota aggiunta nella
    seconda versione, dicembre 2022).}
  \\ Allora dimmi, di che si tratta?
\item[F.] Matematica. 
\item[P.] Uhm, matematica? E da quand'è che ti entusiasma la matematica?
\item[F.] In genere no, ma questa è troppo forte, 
  anche se per me è troppo tardi, visto che le tabelline oramai 
le ho dovute imparare.
\item[P.] Hanno tolto le tabelline dai programmi scolastici?
  Non mi stupirei di niente.
  Tanto oramai l'interesse primario non è quello di formare i ragazzi,
  ma di vendere libri e sussidi, \ldots\  e gite scolastiche, o come
 diavolo le chiamano oggigiorno\ldots \  
 Io alle elementari andavo a scuola
 con il libro di lettura e il sussidiario,
 un quaderno a righe e uno a quadretti. Tu avevi uno zainetto
  che pesava più di te. E alla fine non mi sembra che la
  conoscenza acquisita sia proporzionale al peso dei libri
  e dei quaderni che ci siamo
  portati dietro nella cartella o nello zainetto, 
  o addirittura trascinati nel carrello. 
\item[F.] Non ancora abolite, ma sicuramente lo faranno. Ecco, leggi qua:
  ``Tabelline addio, la moltiplicazione giapponese rende i calcoli un gioco.''
  E ancora, ``Addio tabelline, il metodo giapponese rivoluziona i calcoli.''
\item[P.] Perbacco! E come funziona?
\item[F.] Ecco, guarda, ho visto il video e ti faccio vedere come si
  fa 21 per 23. 
\item[P.] 483, giusto?
\item[F.] Con te non c'è gusto. Hai fatto prima dei giapponesi.
  Ma come hai fatto?
\item[P.]  Non ci vuole molto. 20 volte 23 fa 460, a cui devi aggiungere
  una volta 23. I numeri sono facili e non ci sono riporti.
\item[F.] Sì, ma così devi pensare!
\item[P.] Ah, dimenticavo, voi siete quelli che vorreste non pensare,
  tanto c'è qualcuno che pensa per voi. Uhm, francamente non mi sembra
  una grande pensata la vostra. 
\item[F.] No, non è così\ldots
\item[P.] Va be', lasciamo stare le discussioni sui `massimi sistemi' 
e veniamo al problema
  specifico. Allora, cosa fanno i giapponesi invece di pensare?
\item[F.] Guarda, disegno su un foglio due linee parallele per il primo `2'
 e  poi, ben di\-stanziata, un'altra linea parallela che sta per `1'.
  E le faccio inclinate di 45 gradi, poi vedrai perché.
\mbox{}\vspace{-0.3cm}  
\begin{center}
\epsfig{file=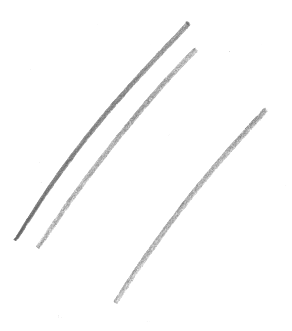,clip=,width=0.25\linewidth}
\end{center}
\mbox{}\vspace{-1.0cm}  
\item[P.] Perché no? Anche se la notazione mi sembra decisamente regressiva
  rispetto a quella escogitata dagli arabi o chi per loro. E direi
  che non è neppure ecosostenibile, come si dice oggi per essere alla moda,
  visto per scrivere
  21 hai consumato un quarto di foglio. Ma vai avanti.
\item[F.] Ora ci traccio sopra
  le linee che indicano 23, ma a 90 gradi rispetto
  a quelle del 21, ecco, così.
\mbox{}\vspace{-0.2cm}  
\begin{center}
\epsfig{file=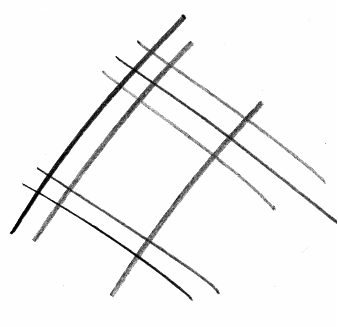,clip=,width=0.3\linewidth}
\end{center}
\mbox{}\vspace{-1.2cm}  
\item[P.] Interessante!
\item[F.] Mi prendi in giro?
\item[P.] No, no, lo sto dicendo sul serio. Ma vai avanti, mi sembra di
  aver capito. 
\item[F.] Certo che hai capito, sicuramente. Io invece
  so come arrivare al risultato ma non ho ancora capito bene perché.
  Ora ti faccio vedere. 
\item[P.] In effetti il risultato è già lì, in forma grafica,
  e va solo letto opportunamente.
\item[F.] Ecco, vediamo le linee dei `2'\ldots
\item[P.] Quelle delle decine\ldots 
\item[F.] \ldots che si intersecano in quattro punti,
  e scriviamo sotto `4'.
\mbox{}\vspace{-0.2cm}  
\begin{center}
\epsfig{file=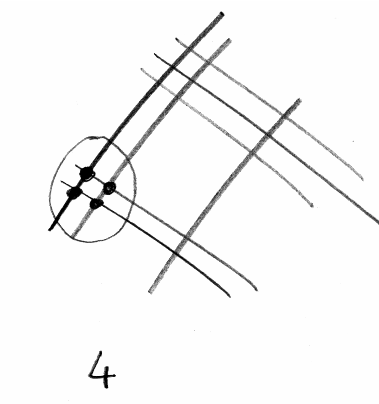,clip=,width=0.35\linewidth}
\end{center}
\mbox{}\vspace{-0.6cm} \mbox{} \\   
  Poi ci sono le linee che indicano `1' e `3'\ldots
\item[P.] Quelle delle unità\ldots 
\item[F.] \ldots che hanno solo tre punti di intersezione, e scriviamo `3'.
\mbox{}\vspace{-0.2cm}  
\begin{center}
\epsfig{file=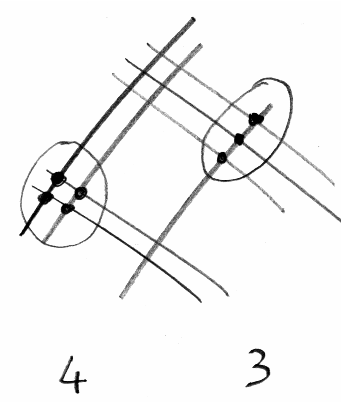,clip=,width=0.35\linewidth}
\end{center}
\mbox{}\vspace{-0.7cm}  \mbox{}\\ 
  E infine le altre intersezioni: le due linee della prima decina -- giusto? --
  che incontrano le tre linee delle unità di `23', e
  le due linee della seconda decina con la linea della prima unità.
  In totale fanno otto intersezioni
\mbox{}\vspace{-0.2cm}  
\begin{center}
\epsfig{file=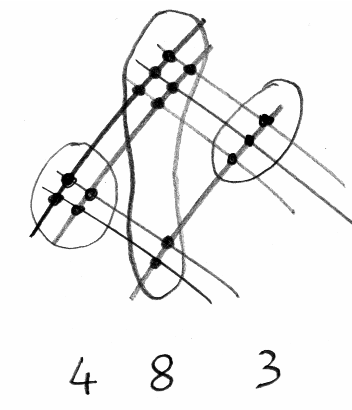,clip=,width=0.35\linewidth}
\end{center}
\mbox{}\vspace{-1.2cm}  \mbox{}\\ 
\item[P.] Simpatico, anche se non direi efficiente, visto il tempo che
  ci hai messo e le risorse che hai dovuto impiegare, intendo carta e
  inchiostro. 
\item[F.] Esagerato! Ma la cosa importante è che non ho dovuto pensare\ldots\ 
  insomma solo un po' per ricordare la regola, 
  e, soprattutto,  non mi è servito sapere le tabelline. 
\item[P.] In effetti! Ma questa moltiplicazione era talmente facile
  che non fa testo. E se le cose si complicano? Per esempio
  se hai numeri a più di due cifre, o\ldots
\item[F.] Ecco, ora ti faccio vedere l'altra che c'era sul video.
  Si tratta di due numeri a 3 cifre, 123 per 321. 
\item[P.] Già, fammi vedere.  
\item[F.] Di nuovo, una linea, poi due linee e poi ancora tre linee.
  Ecco, così:
 \mbox{}\vspace{-0.2cm}  
\begin{center}
\epsfig{file=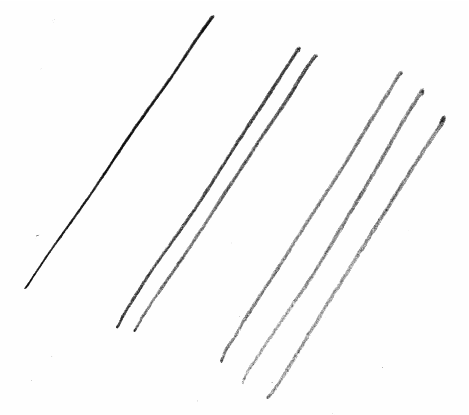,clip=,width=0.35\linewidth}
\end{center}
\mbox{}\vspace{-0.7cm}  \mbox{}\\  
  E abbiamo scritto il primo numero. \\
  Poi, a 90 gradi rispetto alle prime, tre linee, seguite da due e poi da una:
   \mbox{}\vspace{-0.2cm}  
\begin{center}
\epsfig{file=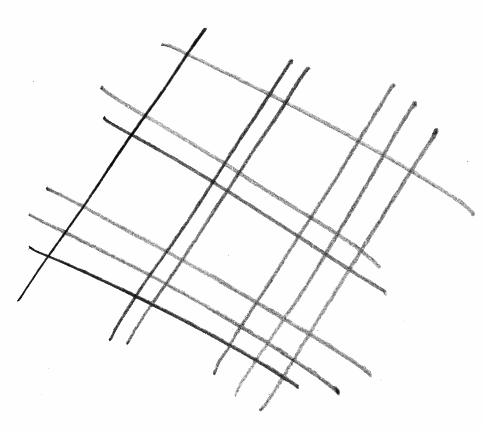,clip=,width=0.35\linewidth}
\end{center}
\mbox{}\vspace{-0.7cm}  \mbox{}\\  
  E ora ci siamo. Basta solo contare il numero di intersezioni.
\item[P.] Buon divertimento! 
\item[F.] In effetti è un divertimento!
\item[P.] Sarà, ma non se devi fare dei conti in tempi rapidi,
  per qualcosa che ti serve.\\
  Ma vai avanti, ti voglio vedere all'opera,
  visto che il gioco ti diverte.
\item[F.] Per prima cosa contiamo le intersezioni sulla destra, che sono tre.
 \mbox{}\vspace{-0.2cm}  
\begin{center}
\epsfig{file=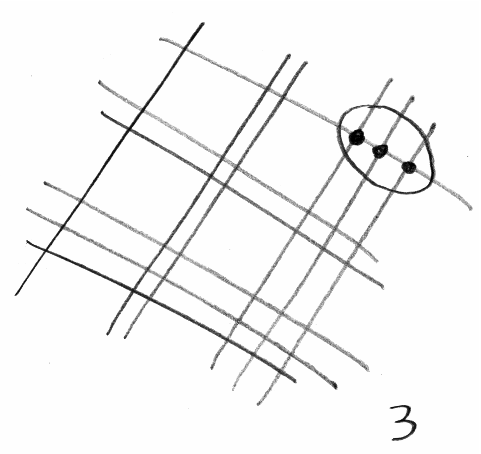,clip=,width=0.35\linewidth}
\end{center}
\mbox{}\vspace{-1.6cm}  \mbox{}\\    
\item[P.] Ok. 
\item[F.] Poi sommiamo tutte le intersezioni allineate a sinistra
  di queste tre, \ldots\ 
  uhm, aspetta, mi sto confondendo, fammi rifare bene la figura, perché
  mi sono persa gli allineamenti.
\item[P.] Ah, andiamo bene\ldots\  e la moltiplicazione è ancora facile\ldots
\item[P.] Dai, papà, sei tu che mi fai emozionare.
  Aspetta, vado a prendere il righello.  
  Ecco, ora si capisce meglio.
 \mbox{}\vspace{-0.2cm}  
\begin{center}
\epsfig{file=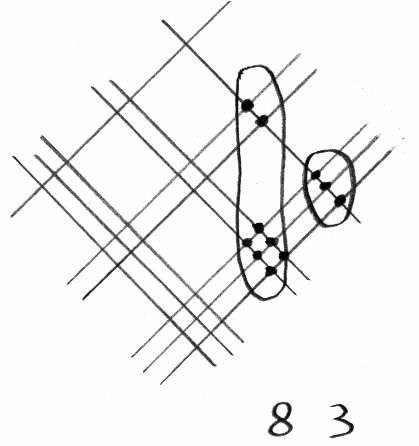,clip=,width=0.35\linewidth}
\end{center}
\mbox{}\vspace{-0.7cm}  \mbox{}\\    
Ho due intersezioni sopra e sei sotto, per un totale di otto.\\
  Poi ho, andando ancora a sinistra, un'intersezione in alto, quattro
  al centro e nove in basso, per un totale di 14.
  \mbox{}\vspace{-0.2cm}  
\begin{center}
\epsfig{file=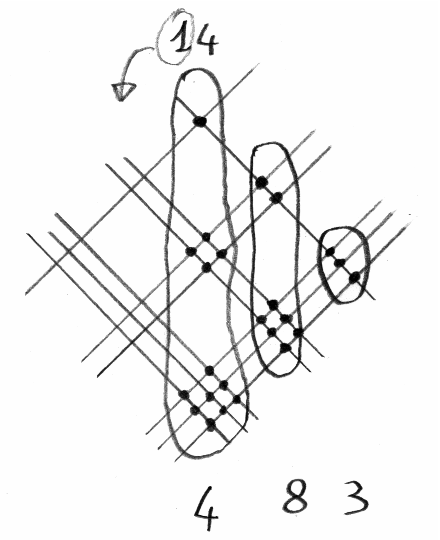,clip=,width=0.35\linewidth}
\end{center}
\mbox{}\vspace{-1.3cm}  \mbox{}\\
\item[P.] Scusa, come hai fatto a capire che le intersezioni
  al centro sono quattro?
\item[F.] Si vedono!
\item[P.] E che quelle in basso sono nove? 
\item[F.] Beh, tre per tre nove!
\item[P.] E questo è il metodo che evitava l'uso
  delle tabelline?
\item[F.] In effetti avrei potuto contare i puntini. 
\item[P.] Giusto! Anche quando calcoli che quattro volte cinque fa venti
  puoi contare, portando alle labbra in successione la punta di ciascun
dito della mano sinistra, usando le dita della mano destra
  per ricordarti quante volte lo hai fatto. Ma vai avanti.
\item[F.] Quindi, dicevo, 14 intersezioni.
  Scrivo allora `4' e mi ricordo di `1',
che devo sommare alle intersezioni a sinistra di queste.
\mbox{}\vspace{-0.2cm}  
\begin{center}
\epsfig{file=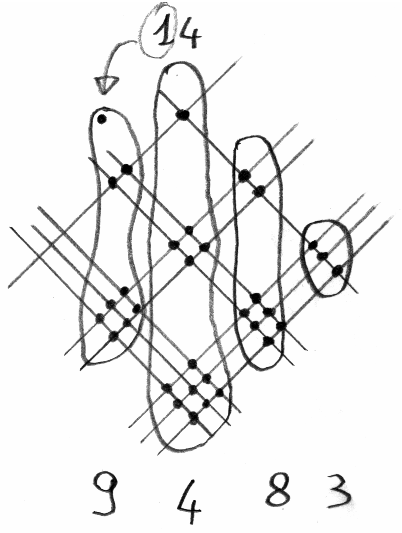,clip=,width=0.35\linewidth}
\end{center}
\mbox{}\vspace{-0.7cm}  \mbox{}\\    
   Ecco, ce ne sono
  due sopra e sei sotto e, sommando `1' che avanzava prima,
  otteniamo `9'.
\item[P.] Brava!  
\item[F.] E infine rimangono le tre intersezioni a sinistra,
\mbox{}\vspace{-0.1cm}  
\begin{center}
\epsfig{file=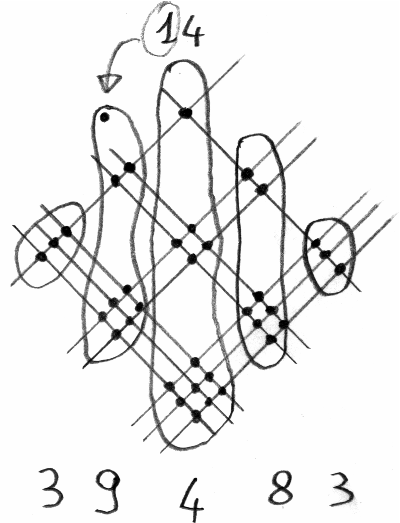,clip=,width=0.35\linewidth}
\end{center}
\mbox{}\vspace{-0.7cm}  \mbox{}\\  
  quindi finalmente\ldots
\item[P.] Già, finalmente\ldots
\item[F.] Non essere polemico, papà. Insomma, abbiamo ottenuto,  
  nell'ordine `3', `9', `4', `8' e `3',
  ovvero {\bf 39483}. 
\item[P.] Il risultato è corretto, bene. Queste erano le moltiplicazioni
  del video, giusto?
\item[F.] Sì, e ti confesso che prima di fartele vedere mi ero già
  allenata un paio di volte. 
\item[P.] Te ne posso proporre una io?
\item[F.] See, immagino! Conoscendoti, sarà come minimo 345679 per 874005,
  o qualcosa del genere\ldots
\item[P.] Non esageriamo, facciamo semplicemente 67 per 85.
  La sapresti fare anche col metodo tradizionale, giusto?
\item[F.] Ci mancherebbe altro! Alle elementari, con la maestra Stefania,
  facevamo cose ben più complicate. 
\item[P.] Ma sapevi già le tabelline. Ora prova a farla con il metodo dei giapponesi
  -- hai detto che si chiama così, giusto? -- e soprattutto senza usare
  le tabelline. 
\item[F.] Facile, sei linee per il `6', \ldots , ma stavolta le
  faccio subito con il righello\ldots \\
  Eccole tutte
  \mbox{}\vspace{-0.2cm}  
\begin{center}
\epsfig{file=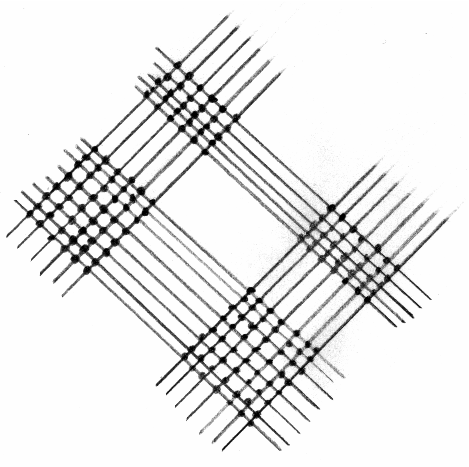,clip=,width=0.4\linewidth}
\end{center}
\mbox{}\vspace{-1.1cm}  \mbox{}\\  
\item[P.] Bello! Bel motivo per una tovaglia, o anche per una sciarpa, questione di gusti. Risultato? 
\item[F.] Allora, cominciamo\ldots\  innanzitutto 35 intersezioni a destra.
\item[P.] Come fai a saperlo? 
\item[F.] Sette per cinque trentacinque.
\item[P.] Avevi detto senza tabelline! Conta! 
\item[F.] See, mi sembra la barzelletta ``zitto e nuota!''\footnote{Famosa
  barzelletta della serie ``zitto e\ldots''. In questo caso è
  la risposta del padre al figlio che piagnucolando dice
``ma papà, io in America non ci voglio andare!''}
\item[P.] Ti ci sei messa tu nei pasticci, veramente ora zitta e conta!
\item[F.] Ma contare zitta non ci riesco!
\item[P.] Va be', basta che conti.
\item[F.] Uno, due, tre, \ldots 
\item[P.] A due a due, niente, eh? 
\item[F.] Va be', ricomincio. Due, quattro, sei,\ldots , 
trentadue, trentaquattro, 
e uno trenta\-cinque. E fortuna che già
  sapevo quanto doveva fare, se no vatti a fidare\ldots\\
  Quindi scrivo `5' e mi segno il `3' da riportare a sinistra.  
\item[P.] Segna, segna e non ti sbagliare.
\item[F.] Ora ci sono le intersezioni centrali. Ma sono un'infinità!
\item[P.] Non esageriamo! 
\item[F.] Sono sei per cinque sopra e otto per sette sotto. 
\item[P.] Non ragionare come se stessi già usando le tabelline. 
          Conta, conta, visto che te la sei cercata!
\item[F.] Papà, facciamo una cosa, facciamo finta che abbia contato,
  ma uso le tabelline. Te ne prego! Se no ci metto una vita e devo ricontare
  i punti almeno due volte per essere sicura di non essermi sbagliata. 
\item[P.] Va bene, ma solo perché sei la mia figlia prediletta\ldots
\item[F.] \ldots non essendocene altre\ldots \ 
  Allora 30 sopra e 56 sotto, che fanno 86 -- fin qui ci arrivo,
  e la maestra Stefania sarebbe fiera di me. \\
  Poi mi devo ricordare dei 3 che avanzavano prima, e siamo arrivati
  a 89, quindi -- finalmeeente! -- scrivo `9' e mi segno `8'. \\
  E per ultimo, abbiamo otto per sei, che fa 48 e che,
  con `8' che riportavo, arrivo a 56. Insomma, se non mi sono sbagliata
  (e per fortuna mi hai risparmiato il supplizio di contare!),
  il prodotto vale 5695. 
\item[P.] Giusto, come posso verificare usando R\,\cite{R} sul mio pc. 
\item[F.] Non c'è gusto! Dovevi verificare con il metodo `nostro'.
\item[P.] Non sono mica scemo! Una cosa è saperlo usare, e capire su cosa
  si basa, e un'altra è usarlo sempre per pedanteria o masochismo.
  Allora, che te ne
  pare del metodo giapponese che rivoluzionerà l'insegnamento della
  matematica alle elementari e eviterà alle generazioni future di imparare
  le tabelline?
\item[F.] Alla fin fine direi che era meglio come avevo imparato
  con la maestra Stefania, usando poi in pratica la calcolatrice
  o quei programmi che usi tu al computer\ldots 
\item[P.] \ldots ma cercando di far girare le nostre rotelline quando
  si tratta di conti con numeri facili, come 20 per 30, o anche
  21 per 23. 
\item[F.] Quindi altra bufala che gira su internet, come tante altre?
\item[P.] Se viene presentata con i titoli che mi hai mostrato, direi
  sicuramente di sì. Se invece è un modo per vedere le moltiplicazioni sotto
  un altro punto di vista, allora la cosa può essere divertente, e anche
  istruttiva.
\item[F.] Cosa vuoi dire?
\item[P.] Ripartiamo dal tuo 21 per 23. Il primo numero è composto da
  due decine e da una unità, il secondo da due decine e tre unità.
  Indicando con $u$ le unità e con $d$ le decine, li possiamo scrivere
  quindi come $(2d + 1u)$  e $(2d + 3u)$. Quando li moltiplichiamo\ldots
\item[F.] Prodotto di binomi\ldots 
\item[P.] Brava, vai avanti tu\ldots 
\item[F.] ``Che a me vien da ridere'', come dici sempre.
\item[P.] È proprio il caso\ldots
\item[F.] Quindi, nell'ordine, abbiamo, $4 d^2$, $6 du$, $2ud$ e
  $3u^2$. Siccome le decine al quadrato sono centinaia, le decine per le unità
  sono decine, indipendentemente dall'ordine, e le unità al quadrato
  sono delle unità, otteniamo 4 centinaia, 8 decine e 3 unità,
  insomma 483. Quindi -- vediamo se ho capito -- nel caso
  di moltiplicazioni a due cifre, il metodo giapponese dà una rappresentazione
  grafica del prodotto di due binomi. E così via per moltiplicazioni a più
  cifre, giusto?
\item[P.] Giustissimo! E, a proposito, riguardando i tuoi disegni
  con le linee, quanti raggruppamenti c'erano quando facevi 21 per 23 e
  67 per 85?
\item[F.] Quattro: quello a destra, quello a sinistra e i due centrali. 
  Proprio i quattro termini del prodotto di due binomi.
\item[P.] E quando facevi 123 per 321?
\item[F.] Aspetta che riguar\ldots
\item[P,] No, senza guardare, pensa un attimo. 
\item[F.] Uno a destra, uno a sinistra, due\ldots 
\item[P.] Ti ho chiesto di pensare, non di ricordare\ldots\  Anzi
  ti complico il problema, con la speranza di facilitarti la visione
  d'insieme. Se avessi dovuto fare 12345678 per
 87654321, insomma il prodotto di due numeri a otto cifre, 
 quanti raggruppamenti avresti avuto?
\item[F.] Aspe'\ldots
\item[P.] Lascia stare righello e matita, e pensa a una scacchiera. 
\item[F.] Cambiamo gioco?
\item[P.] No, analizziamo un  problema analogo, caratterizzato
  da una dimensionalità uguale a quella che ci interessa.
  Insomma, eravamo partiti
  dal prodotto di due numeri a due cifre e poi siamo passati
  a quello di due numeri a tre cifre. Ora siamo saltati a otto. 
\item[F.] Ah, ecco, dove vuoi arrivare! La scacchiera ha otto righe e
  otto colonne, per un totale di sessantaquattro caselle,
  insomma sessantaquattro possibilità, almeno per la regina.
  Rappresentano quindi 
  tutti gli incroci possibili -- papà, li posso chiamare incroci? --
  fra le righe e le colonne.
\item[P.] Chiamali come ti pare, se hai capito di cosa si tratta. 
 \item[F.] Quindi -- ecco! -- il prodotto di due numeri da otto cifre l'uno 
  produce, con il metodo dei giapponesi, sessantaquattro incroci
  di righe e colonne, cioè sessantaquattro raggruppamenti. Mentre
  nel caso di numeri a tre cifre ne producono nove, tre per tre.
  Sì, erano nove, ora che ci penso bene.
\item[P.] E anche con il nostro metodo abbiamo lo stesso 
  numero di `raggruppamenti', per così dire. 
\item[F.] Ma nel nostro metodo non ci sono raggruppamenti!
\item[P.] Ah, no? Puoi fare la moltiplicazione `123 per 321'
  come avevi imparato alle elementari?
\item[F.] Un momento. Allora,  uno per tre, uno per due,\ldots;
 due per tre, due per due, \ldots \  Ecco, 39483.
\item[P.] E ora metti i due conti uno affianco all'altro, 
  intendo quello nipponico e quello {\em de noantri}. 
\item[F.] Eccoli.
  \mbox{}\vspace{-0.2cm}  
\begin{center}
  \epsfig{file=Fig12.eps,clip=,width=0.35\linewidth}
  \hspace{1.3cm}
  \epsfig{file=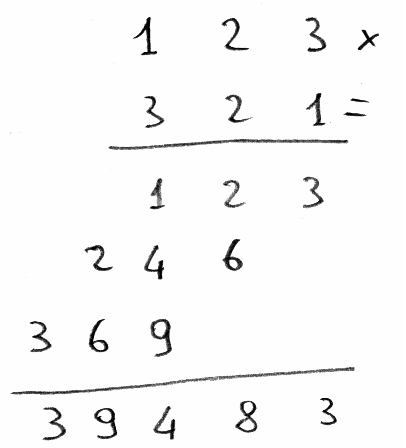,clip=,width=0.35\linewidth} 
\end{center}
\mbox{}\vspace{-0.3cm}  \mbox{}\\ 
  Ma è la stessa cosa! Il `3' a destra da solo, poi i `2' e il `6' 
  allineati, che fanno otto, e così via. 
\item[P.] Toh, curioso, non ti pare? Ma forse non sarà un caso\ldots
\item[F.] Ma in `67 per 85' non funziona.
\item[P.] E perché no? La differenza è che quando il numero
  di intersezioni in ciascun raggruppamento  
  supera 9, noi facciamo i riporti al volo. Prova a fare 
  67 per 85 con il nostro metodo senza riporti.
\item[F.] Come si fa? Viene un casino. 
\item[P.] Basta che metti i risultati del singolo prodotto 
  fra parentesi. 
\item[F.] Ecco, va bene così?
   \mbox{}\vspace{-0.2cm}  
\begin{center}
  \epsfig{file=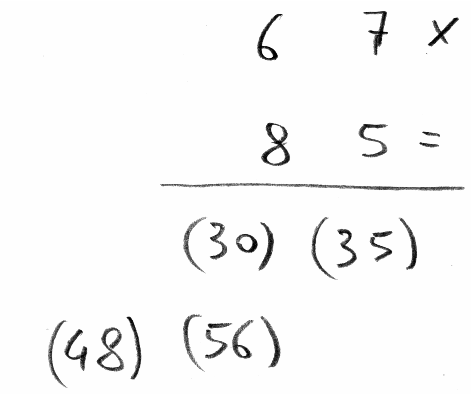,clip=,width=0.35\linewidth} 
\end{center}
\mbox{}\vspace{-1.4cm}  \mbox{}\\ 
\item[P.] Va bene per forza, se non hai sbagliato a fare i conti.
\item[F.] Papà, le tabelline le so ancora! 
\item[P.] Bene, ora prova a rivedere il risultato 
  alla luce di quanto abbiamo detto prima. 
\item[F.] In effetti, abbiamo di nuovo quattro raggruppamenti. 
  Uno da 35 unità, che è quello a destra in alto,
  uno da 48 centinaia, che è quello a sinistra in basso
  e due, rispettivamente da 30 e da 56 decine, che sono i
  due incolonnati. 
\item[P.] Ergo?
\item[F.] Ergo sum\ldots\  lasciamo stare cosa\ldots \ 
  Ma stavolta ho finalmente capito!
  In effetti 35 unità sono cinque unità e tre decine. E scrivo `5'. 
  Poi abbiamo 86 più 3 decine, per un totale di 89 decine, 
  ovvero 9 decine e 8 centinaia. Quindi scrivo `9' a sinistra di `5'
  e vado avanti. Infine abbiamo 48 più 8 centinaia, ossia 56 centinaia, 
  ovvero sei centinaia e cinque migliaia. Insomma {\bf 5695}. 
  Stesso risultato, e ho pure capito perché e cosa significano 
  le singole cifre, perché, a ripensarci, quel metodo 
  giapponese, come dicono, oltre a essere praticamente inutilizzabile 
  se le cifre sono grandi, è veramente da automi. 
\item[P.] E il numero di operazioni da fare? 
\item[F.] Lo stesso, perché devi fare una moltiplicazione
  per ogni raggruppamento, e il numero dei raggruppamenti
  è pari al prodotto dei numeri di cifre di ciascun fattore.
  Infine devi fare le somme dei raggruppamenti delle stesse `cose',
  insomma delle decine, delle centinaia, e così via. E poi ci sono
  i riporti. Insomma i conti sono gli stessi. 
\item[P.] E con le tabelline come la mettiamo?
\item[F.] Mi sembrano inevitabili, a meno di non usare 
  solo calcolatrici o programmi al computer. Fare a mano
  il conteggio del numero di intersezioni in ciascun raggruppamento
  mi sembra una follia, senza contare i disegni, gli allineamenti 
  e la possibilità di sbagliarsi.  
  Ma ora basta perché è dai tempi della maestra Stefania
  che non facevo tutti questi conti e mi sta scoppiando la testa. 
\item[P.] E già, perché dopo le elementari, viene la 
  `Terra di mezzo', che sono le medie, e poi la `Grande ubriacatura'
  dei licei --- mi riferisco ai programmi di matematica --- 
  con cose assolutamente inutili o accessorie alle quali
  si dà la stessa importanza delle poche cose fondamentali che 
  dovrebbero rimanere impresse nella mente.\footnote{Per altre
    critiche dell'autore ai programmi scolastici delle superiori,
  si veda Ref.\cite{StregheBayes}.}
  Ma lasciamo stare\ldots
\item[F.] Sì, sì, lasciamo stare\ldots
\item[P.] Aspe', non ho finito. Ogni tanto una paternale 
  ci vuole, è mio diritto e dovere, no? 
\item[F.] Giusto, e io ascolto compìta. Mio dovere. 
\item[P.] Volevo dire che mi ha fatto piacere che, giocandoci
  insieme su questa cosa, ci hai riflettuto un po', 
  e poi mi sembra che ti sia pure divertita.
  Perché la paura mia è nei confronti di coloro che non vogliono che
  pensiate, tanto pensano loro, e non di certo ai vostri interessi.
  O, in tempi di internet, di quelli che dicono cose senza
  averci pensato, e pretendono che siccome ``siamo in democrazia'' la loro
  opinione valga come quella di chi, coscienziosamente, sulle cose
  ci ha meditato non poco, e che in fondo, proprio perché
  intellettualmente onesto, qualche dubbio ce lo può
  ancora avere. Ignoranza e ostentazione di sicurezza sembra che vadano
  troppo spesso insieme, di questi tempi.
  E, a proposito --- e poi ti lascio a quello che stavi facendo ---
  permettimi di leggerti una cosa che mi è arrivata poco fa,
  per mail, da un amico, pescata non so dove nel mare magnum
  dei social. Non sono di quelli che detesta o criminalizza questi
  nuovi mezzi di comunicazione,
  ma ci vuole un po' di\ldots, anzi tanta, accortezza per non affogarcisi,
  o bere schifezze, in questo mare. Un momento che la ritrovo. 
  Ecco, ``il problema della libertà
  di pensiero non è l'assenza di libertà, ma l'assenza di PENSIERO.''
\item[F.] Bella! Però, scusa papà, ma a pensarci bene
  in tante parti del mondo non è proprio così\ldots
\item[P.] E già\ldots\ hai perfettamente ragione. In effetti,
  scusa, ti ho riportato affrettatamente e in modo acritico
  un'affermazione che suona universale, perché l'avevo istintivamente
  riferita alla realtà delle nostre parti. Vedi che bisogna
  stare ve\-ramente attenti?
\end{description}

\begin{center}
  \mbox{}\\
  ------------------------------ 
  \end{center}

 \mbox{}\\  \mbox{}\\

\end{document}